\documentclass[11pt]{article}
\usepackage{amsfonts}
\newcommand{\beq}{\begin{equation}}
\newcommand{\enq}{\end{equation}}

\newtheorem{Theorem}{Theorem}[section]
\newtheorem{Lemma}[Theorem]{Lemma}

\newtheorem{Remark}[Theorem]{Remark}

\newcommand{\benu}{\begin{enumerate}}
\newcommand{\beqa}{\begin{eqnarray}}
\newcommand{\beqan}{\begin{eqnarray*}}
\newcommand{\eay}{\end{array}}
\newcommand{\edm}{\end{displaymath}}
\newcommand{\eenu}{\end{enumerate}}
\newcommand{\eeq}{\end{equation}}
\newcommand{\eeqa}{\end{eqnarray}}
\newcommand{\eeqan}{\end{eqnarray*}}

\newcommand{\br}{\begin{Remark}}
\newcommand{\er}{\end{Remark}}

\newcommand{\bqa}{\begin{eqnarray}}
\newcommand{\eqa}{\end{eqnarray}}
\newcommand{\bqw}{\begin{eqnarray*}}
\newcommand{\eqw}{\end{eqnarray*}}

\newcommand{\non}{\nonumber}
\newcommand{\bea}{\begin{array}{cc}}
\newcommand{\ena}{\end{array}}

\RequirePackage{amsfonts}
\usepackage{stmaryrd}
\usepackage{amssymb}
\usepackage{mathrsfs}
\usepackage{amsmath}
\usepackage{subeqnarray}
\usepackage{verbatim}
\usepackage{cases}
\usepackage{amsfonts}
\usepackage{amssymb}
\usepackage{amssymb}
\usepackage{amssymb}
\usepackage{hyperref}
\begin{document}

\newpage
\pagenumbering{arabic} \setcounter{page}{1}

\begin{center}
	{\Large  Global attractor for a nonlinear one-dimensional compressible viscous micropolar fluid}\\\vspace{0.25in}Lan Huang\footnote[3]{College of Mathematics and Statistics,
		North China University of Water Resources and Electric Power,
		Zhengzhou 450011, Henan,  P. R. China. Email: huanglan82@hotmail.com}\ \ \ \
	\ \ \  \  Xin-Guang Yang   \footnote[4]{Department of Mathematics and Information Science, Henan Normal University, Xinxiang, 453007, Henan,
		P. R. China. Email: yangxinguang@hotmail.com}
	  \ \ \ \ \
	  Yongjin Lu \footnote[2] {Department of Mathematics and Economics, Virginia State University, Petersburg, VA 23806, USA. Email: ylu@vsu.edu}\ \ \ \
	  Taige Wang \footnote[1]{Department of Mathematical Sciences, University of Cincinnati, Cincinnati, OH 45221, USA. Corresponding to: wang2te@ucmail.uc.edu}
	\ \ \ \ \vspace{0.06in}
\end{center}
\vspace{0.06in}
\begin{abstract}
	This paper considers the dynamical behavior of solutions of constitutive systems for 1D compressible viscous and heat-conducting micropolar fluids. With proper constraints on initial data, we prove the existence of global attractors in generalized Sobolev spaces $H_{\delta}^{(1)}$ and $H_{\delta}^{(2)}$.
	
\end{abstract}

\section{Introduction}\label{sec1}
 \ \ \ \ \ The microfluid model was developed by Eringen (see, e.g., \cite{e1, E}) in 1960s. The model describes microscopic phenomena of materials possessing microstructures. The particles in small volume elements have micromotions, for instance, microrotations. As the stress and body moments are coupled with the spin inertia, this type of the constitutive system turns out to be complicated to mathematical analysis. A ``simplified" class of this type of models is the micropolar fluid model, in which the first stress moments and gyration tensor are skew symmetric. The micropolar fluids include a class of anisotropic polymeric fluids which have dumbbell molecules, such as liquid crystals, and blood. For similar types of complex fluids, we refer readers to review monograph \cite{Re}.   \\

When solvents (Newtonian fluids) contain low concentration of polymeric additives, Navier-Stokes constitutive laws can be coupled with microplar constitutive relation to describe their fluid dynamics. In general, the compressible micropolar fluid models can be obtained from integral form of conservation laws, which are coupled with various constitutive relations, such as Fourier's law, Boyle's law and polytropy (See, e.g., N. Mujakovi\'{c} \cite{M1}). \\

The main content of this paper focuses on existence of global attractors for the compressible viscous and heat-conducting micropolar fluid, in a thermodynamical sense: perfect and polytropic. The model is written in terms of Lagrangian coordinate: \begin{equation}\label{1-1}
\begin{cases}
u_t=v_x,& \\
v_t=\left(-p+
\frac{v_x}{u}\right)_x,& \\
\omega_t=A\left((\frac{\omega_x}{u})_x-u\omega\right),& \\
C_v\theta_t=-K\frac{\theta v_x}{u}+\frac{v_x^2}{u}+\frac{\omega^2_x}{u}+u\omega^2+D(\frac{\theta_x}{u})_x,& \end{cases}
\end{equation}
where $t>0$ is time and  $x\in\Omega=[0,1]$ denotes the mass variable. Here
the unknown vector  $(u(x,t),v(x,t), \omega(x,t),\theta(x,t))$ represents
the specific volume ($u=\frac{1}{\rho})$, velocity,  microrotation
velocity, and the absolute temperature of the fluid flow respectively;  the pressure
$p=K\theta/u$ and $K,\;D,\;A,\;C_v$ are positive constants.\\

We consider system \eqref{1-1} subject to the following boundary
condition \begin{eqnarray}  v(0,t)=v(1,t)=0,\; \omega(0,t)=\omega(1,t)=0,\;
\theta_x(0,t)=\theta_x(1,t)=0,\quad t>0\label{1-2}\end{eqnarray} and the initial condition \begin{eqnarray}
(u,v,\omega,\theta)(x,0)=(u_0,v_0,\omega_0,\theta_0)(x) \label{1-3}\end{eqnarray} for
$x\in\Omega=[0,1]$, where $u_0=\frac{1}{\rho_0},\;v_0,\;\omega_0$
and $\theta_0$ are prescribed functions.\\

Moreover, we assume the compatibility condition
\begin{eqnarray}\label{1-4} v_0(0)=v_0(1)=0,\quad\omega_0(0)=\omega_0(1)=0,\quad\theta_{0x}(0)=\theta_{0x}(1)=0\end{eqnarray}
holds.\\

There are plenty of works on existence of solutions of different types of initial-boundary value problems for incompressible fluids (see, e.g., \cite{bl96, L}), but existence theory for compressible micropolar fluid is still in development.\\

The well-posedness result of system \eqref{1-1} is summarized as follows:\\

(1) Under suitably prescribed initial data for 1D micropolar fluid model, N. Mujavoki\'{c} \cite{M1} established the global existence and asymptotic behavior of the solution for the system \eqref{1-1} with the boundary conditions \eqref{1-2} in \cite{M2,M3}, then the authors obtained the exponential stability in \cite{HN, M5,M6} and established the local existence and global existence for the same system with non-homogeneous boundary conditions for velocity and microrotation:

 \begin{eqnarray}\label{1-5}
&&v(0,t)=\mu_0(t),\;v(1,t)=\mu_1(t),\;\omega(0,t)=\nu_0(t),\nonumber\\
&&\omega(1,t)=\nu_1(t),\theta_x(0,t)=\theta_x(1,t)=0.\end{eqnarray}

Recently,
Mujakovi\'{c} in \cite{M7,M8} and references cited therein
studied  the local and global existence  for the system (1.1)-(1.4) with a non-homogeneous boundary condition for temperature  \begin{eqnarray}\label{1-6}
&&v(0,t)=v(1,t)=0,\;\omega(0,t)=
\omega(0,t)=0,\nonumber\\
&&\theta_x(0,t)=\mu_0(t),\theta_x(1,t)=\mu_1(t).\end{eqnarray}
Mujakovi\'{c} and  \v{C}rnjari\'{c}-\v{Z}ic \cite{MC2} proved the global existence with the boundary condition \begin{eqnarray}\label{1-7}
&&v(0,t)=0,\;\omega(0,t)=\omega(1,t)=0,\nonumber\\
&&\theta_x(0,t)=\theta_x(1,t)=0,(\frac{v_x}{u}-K\frac{\theta}{u})(1,t)=0.\end{eqnarray}

(2) In three dimensional domains, for the spherically symmetric motions of compressible micropolar fluids in bounded annular domains, \cite{qiu3, HK,HL} obtained the global existence, the uniqueness, and  asymptotic behavior, the exponential stability, and regularity of generalized solutions.  For the cylindrically symmetric motion in the bounded subset domain of $\mathbb{R}^3$ with two coaxial
cylinders that present the solid thermoinsulated walls, one can refer to  \cite{dr2,dr1, MC1,qw} for similar results.\\

The theory of infinite dimensional dynamical systems involves Navier-Stokes equations, MHD systems, Boussinesq equations, etc.. Since 1980s, long time behaviors of solutions, such as existence of attractors and their geometric structures, are investigated to approach chaos and sigularites present in turbulence. In the past decades, there are plentiful literatures to  deal with this dynamics of 2D incompressible micropolar flows (or its extended models such as magnet-micropolar fluids):

\begin{eqnarray}\label{1-1a}
\begin{cases}
v_t-(\nu+\kappa)\Delta u-2\kappa \nabla\times w+\nabla p+v\cdot\nabla v=f(t),&\\
w_t-\gamma\Delta w+4\kappa w-2\kappa\nabla\times v+v\cdot w=g(t),&\\
\nabla\cdot v=0. &\end{cases}
\end{eqnarray} The forward and pullback attractors and their structures and dimensions for system \eqref{1-1a} in smooth or non-smooth domains can be found in \cite{bl, ch, ccd06, ccd07, cdc07, dc, lu01, lu03, lu04, lu05, mp, ta, ta10, zs, zsh}. If Mach number is close to $0.3$, the incompressible and compressible fluid are almost identical; however, Mach numbers of most of fluids are not around 0.3, and hence there exist the huge differences between compressible and incompressible fluids. So as we know, there are no results available on the existence of attractors for compressible models, even in one dimension.\\

The objective of this paper is to investigate the long-time dynamics of problem \eqref{1-1}-\eqref{1-3} by using the abstract analysis technique established in \cite{6, qinshu, teman,zq1}.
The main features and difficulties of this paper are stated as follows:\\

(I) The first difficulty is to obtain the attracting property via absorbing set in appropriate metric spaces.
In order to prove the existence of absorbing set, we must prove that the orbit of solutions starting from any bounded set of closed subspace will re-enter this closed subspace and stay there after a finite time, which should be uniform for all solution orbits starting from there.\\

Inspired by \cite{qinshu}, we introduce spaces \begin{eqnarray} H^{(1)}&=&\Big\{(u,v,\omega,\theta)\in H^1[0,1]\times H^1[0,1]\times H^1[0,1]\times H^1[0,1]: \nonumber\\ &&u(x)>0,\theta(x)>0,x\in[0,1],
v(0,t)=v(1,t)=0, \omega(0,t)=\omega(1,t)=0\Big\},\non\end{eqnarray}
and \begin{eqnarray}
H^{(2)}&=&\Big\{(u,v,\omega,\theta)\in H^2[0,1]\times H^2[0,1]\times H^2[0,1]\times H^2[0,1]:u(x)>0,\theta(x)>0,x\in[0,1],\nonumber\\ &&
v(0,t)=v(1,t)=0, \omega(0,t)=\omega(1,t)=0,\theta_x(0,t)=\theta_x(1,t)=0\Big\}\non\end{eqnarray}
which become two metric spaces equipped with the metrics induced from the usual norms.
$H^1$ and $H^2$ are the usual Sobolev spaces in the above.
Let $\delta_i(i=1,2,3,4,5)$ be any given constants such that \begin{eqnarray}\label{2-1}
\delta_1\in R,0<\delta_5<\delta_2,\delta_3>0,
\delta_4\geq\max\left[\frac{{e}^{\frac{\delta_1}{K}}}{2(2\delta_2/{C_V})^{\frac{C_V}{K}}},\delta_3\right]>0\end{eqnarray}
and let \begin{eqnarray} H_{\delta}^{(i)}&=&\bigg\{(u,v,\omega,\theta)\in H^{(i)}:\int_0^1(C_V\log\theta+K\log u)dx\geq \delta_1,\non\\
&&\delta_5\leq\int_0^1(C_V\theta+\frac{v^2}{2}+\frac{\omega^2}{2A})dx\leq\delta_2,
\delta_3\leq\int_0^1udx\leq\delta_4,\non\\
&&\frac{\delta_5}{2C_V}\leq\theta\leq\frac{2\delta_2}{C_V},\frac{\delta_3}{2}\leq u\leq 2\delta_4\bigg\},\quad i=1,2.\non\end{eqnarray}
Obviously, $H_{\delta}^{(i)}\;(i=1,2)$ is a sequence of closed subsequences of $H^{(i)}(i=1,2)$.\\

Using some results from \cite{1,HN,M3,teman} and more
precise estimates to deal with more complex terms, we can prove the
existence of $C_0$-semigroup in $H^{(1)}$ for problem \eqref{1-1}-\eqref{1-3} and obtain the absorbing set in $H^{(1)}_{\delta}$.\\

(II)  The second difficulty is that the first three constraints in (I) are invariant (Lemma \ref{le4-1}), while the last two constraints are not invariant (Lemma \ref{le4-2}). Since the original spaces $H^{(i)}(i=1,2)$ are incomplete, we use $H_{\delta}^{(i)}\;(i=1,2)$ introduced in (III) to overcome this obstacle.\\

(III)  By virtue of idea from \cite{qinshu} and delicate uniform estimates, we obtain the existence of global attractors in $H^{(i)}_{\delta} (i=1,2)$ which are compact, invariant $\omega$-limit sets, see Theorem \ref{th2-1}.\\

The rest of this paper is organized as follows: the main result is stated in Section \ref{sec2}; in Section \ref{sec3},
we shall give the proof that operators $\{S(t)\}$ defined by the solutions form a $C_0$-semigroup on $H^{(i)}\;(i=1,2)$; in Sections \ref{sec4} and \ref{sec5}, we shall
establish the existence of an absorbing set in respectively $H^{(1)}_{\delta}$ and $H^{(2)}_{\delta}$,  then finalize the proof of Theorem \ref{th2-1} in Section \ref{sec6}.

\section{Main result}\label{sec2}

The notation in this paper is shown as follows:\\

$L^{\bar{p}}, 1\leq \bar{p}\leq +\infty, W^{m,\bar{p}}, m\in N,
H^1=W^{1,2}, H_0^1=W_0^{1,2}$ denote the usual (Sobolev) spaces on
$(0, 1).$ In addition,  $\|\cdot\|_{B}$ denotes the norm in the
space $B$; we also put $\|\cdot\|= \|\cdot\|_{L^2}.$  Subscripts $t$
and $x$ denote the (partial) derivatives with respect to $t$ and
$x$, respectively. We use $C_0^{(i)}(i=1,2)$ to denote the generic
positive constant depending only on  $H^i$ norm of initial
datum
$(u_0,v_0,\omega_0,\theta_0)$, $\displaystyle{\min_{x\in
		[0,1]}u_0(x)}$ and $\displaystyle{\min_{x\in [0,1]}\theta_0(x)}$,
but independent of variable $t$. $C_{\delta}$ or $C_{\delta}'$ denotes the universal constant depending
only on $\delta_i$'s $(i=1,2,3,4,5)$, but independent of initial data. $C_{\delta}^{(i)}(i=1,2)$ depending on both $\delta_j$'s $(j=1,2,3,4,5)$, the $H^i$ norm of  the initial data $(u_0,v_0,\omega_0,\theta_0)$, $\displaystyle{\min_{x\in
		[0,1]}u_0(x)}$ and $\displaystyle{\min_{x\in [0,1]}\theta_0(x)}$. $C$ denotes the generic absolute positive constant independent of $\delta$ and the initial data.\\

Now we can state our main result as following.
\begin{Theorem}\label{th2-1}
	The nonlinear semigroup $S(t)$ defined by the solution to problem \eqref{1-1}-\eqref{1-3} maps $H^{(i)}\;(i=1,2)$ into
	itself. Moreover, for any $\delta_i(i=1,2,3,4,5)$, it possesses a maximal  universal attractor ${\mathcal{A}}_{i,\delta}$ in $H_{\delta}^{(i)}\;(i=1,2)$.
\end{Theorem}

By the theory of global attractors in \cite{6, qinshu, teman}, we can see that an $\omega$-limit set is a global attractor if it is nonempty, compact, invariant for the continuous semigroup. These property can be achieved by proving (a) continuity of semigroup, (b) compactness, (c) attracting property. In this paper, we shall verify the continuity of semigroup in Section \ref{sec3} and prove the attracting set in Sections \ref{sec4}--\ref{sec5}, and the compactness via the compact embedding of Sobolev spaces and uniform energy estimates.

\section{ Nonlinear $C_0$-semigroup on $H^{(i)}\;(i=1,2)$ }\label{sec3}
\setcounter{equation}{0}
As mentioned in the previous section, for any initial data $(u_0,v_0,\omega_0,\theta_0)\in H^{(1)}$, the results on global existence, uniqueness, and asymptotic behavior of solutions to problem \eqref{1-1}-\eqref{1-3} have been established in \cite{HN,M3}, respectively.

\begin{Lemma}\label{le3-1} Assuming the initial data $(u_0,v_0,\omega_0,\theta_0)\in H^{(1)}$ and the compatibility condition \eqref{1-4} are satisfied,
	then there exists a unique generalized global solution $(u(t),v(t),\omega(t),\theta(t))$ in $H^{(1)}$ to the problem \eqref{1-1}-\eqref{1-3} which satisfies \begin{equation}\label{3-1}
	\begin{cases} 0<1/C_0^{(1)}\leq u(x,t)\leq C_0^{(1)},&\\
	0<1/C_0^{(1)}\leq \theta(x,t)\leq C_0^{(1)}\;\;{\rm on}\;\;\;
	\Omega\times (0,\infty),&\\
	u_x\in L^{\infty}(0,\infty;L^2(\Omega))\cap L^2(0,\infty;L^2(\Omega)),&\\
	 v,\omega \in L^{\infty}(0,\infty;H^1(\Omega))\cap L^2(0,\infty;H^2(\Omega)),&\\
	v_t,\omega_t\in L^2(0,\infty;L^2(\Omega)),&\\
	\theta\in L^{\infty}(0,\infty;H^1(\Omega)),&\\
	\theta_x\in L^{2}(0,\infty;H^1(\Omega)),\;\theta_t\in L^2(0,\infty;L^2(\Omega)),&
	\end{cases}\end{equation}
	and there exist positive constants $\tilde{\gamma_1}=\tilde{\gamma_1}(C_0^{(1)}),C_0^{(1)}$ such that, for any fixed $\tilde{\gamma}\in(0,\tilde{\gamma_1}]$ and for any $t>0$,
	\begin{eqnarray}\label{3-2}&& e^{\gamma t}(\|u-\bar
	u\|^2_{H^1}+\|v\|_{H^1}^2+\|\omega\|_{H^1}^2+\|\theta-\bar
	\theta\|^2_{H^1})\non\\
	&&+\int_0^te^{\gamma s}(\|u-\bar
	u\|^2_{H^1}+\|v\|_{H^2}^2+\|\theta-\bar
	\theta\|^2_{H^2}+\|\omega\|_{H^2}^2+\|v_t\|^2+\|\omega_t\|^2+\|\theta_t\|^2)(s)\,ds\non\\
&&\leq
	C_0^{(1)},\end{eqnarray} where $\bar u= \int_0^1
	u(x)\,dx=\int_0^1 u_0(x)\,dx$, $\bar\theta=\frac{1}{C_v}\int_0^1(\frac{1}{2}v_0^2+\frac{1}{2A}\omega_0^2+C_v\theta_0)\,dx$.
\end{Lemma}
{\bf Proof.}  The estimate \eqref{3-1} and \eqref{3-2} were obtained in \cite{M3} and \cite{HN} respectively. The proof is complete. \hspace*{\fill}$\Box$

\begin{Lemma}\label{le3-2} Assuming initial data $(u_0,v_0,\omega_0,\theta_0)\in H^{(2)}$ and the compatibility condition \eqref{1-4} are satisfied, then there exists a unique generalized global solution $(u(t),v(t),\omega(t),\theta(t))$ in $H^{(2)}$ to the problem \eqref{1-1}-\eqref{1-3} which satisfies \begin{eqnarray}\label{3-3}&&\|u(t)-\bar u\|_{H^2}^2+\|v\|_{H^2}^2+\|\omega\|^2_{H^2}
	+\|\theta-\bar\theta\|_{H^2}^2\nonumber\\
	&&+\int_0^t(\|u(t)-\bar u\|_{H^2}^2+\|v\|_{H^3}^2+\|\omega\|^2_{H^3}
	+\|\theta-\bar\theta\|_{H^3}^2+\|v_t\|_{H^{1}}^2+\|\omega_t\|_{H^{1}}^2+\|\theta\|_{H^{1}}^2)(s)ds\non\\
	&&\leq C_0^{(2)},\end{eqnarray}
	and there exists positive constant $\tilde{\gamma_2}=\tilde{\gamma_2}(C_0^{(2)}),C_0^{(2)}>0$, such that for any fixed $\tilde{\gamma}\in(0,\tilde{\gamma_2}]$ and for any $t>0$,
	\begin{eqnarray}\label{3-4}
	&& e^{\gamma t}(\|u-\bar
	u\|^2_{H^2}+\|v\|_{H^2}^2+\|\omega\|_{H^2}^2+\|\theta-\bar
	\theta\|^2_{H^2})\nonumber\\
	&&+\int_0^te^{\gamma s}(\|u-\bar
	u\|^2_{H^2}+\|v\|_{H^3}^2+\|\theta-\bar
	\theta\|^2_{H^3}+\|\omega\|_{H^3}^2+\|v_t\|_{H^1}^2+\|\omega_t\|_{H^1}^2+\|\theta_t\|_{H^1}^2)(s)\,ds\non\\
	&&\leq
	C_0^{(2)}.\end{eqnarray}
\end{Lemma}
{\bf Proof. }  The estimate \eqref{3-3}-\eqref{3-4} were obtained in \cite{HN}. \hspace*{\fill}$\Box$

\begin{Lemma}\label{le3-3} The unique generalized global solution $(u(t),v(t),\omega(t),\theta(t))$ in $H^{(1)}$
	defines a nonlinear $C_0$-semigroup $S(t)$ on $H^{(1)}$. Moreover, for any $(u_0,v_0,\omega_0,\theta_0)\in H^{(1)}$, the generalized global solution $(u(t),v(t),\omega(t),\theta(t))$ to the problem \eqref{1-1}-\eqref{1-3} satisfies
	\begin{eqnarray}
	&&(u(t),v(t),\omega(t),\theta(t))=S(t)(u_0,v_0,\omega_0,\theta_0)\in C([0,\infty), H^{(1)}),\label{3-5}\\
	&&u(t)\in C^{\frac{1}{2}}([0,+\infty),H^1),\quad v(t),\omega(t),\theta(t)\in C^{\frac{1}{2}}([0,+\infty),L^2).\label{3-6}
	\end{eqnarray}
\end{Lemma}
{\bf Proof. } We will separate the proof into two steps. We will first show the semigroup $\{S(t)\}$ is uniformly bounded in $H^{(1)},$ then prove the continuity of $\{S(t)\}.$\\

\emph{Step 1.} By Lemma \ref{3-1}, we know that for any $t>0$, the operator $S(t):(u_0,v_0,\omega_0,\theta_0)\in H^{(1)}\mapsto (u(t),v(t),\omega(t),\theta(t))\in H^{(1)}$ for any $t>0$ exists and, by the uniqueness of generalized global solutions, satisfies on $H^{(1)}$, for any $t_1,t_2\in[0,\infty)$,
\begin{eqnarray} S(t_1+t_2)=S(t_1)S(t_2)=S(t_2)S(t_1).\label{3-7}\end{eqnarray}
Moreover, by Lemma \ref{3-1}, $S(t)$ is uniformly bounded on $H^{(1)}$ with respect to $t>0$, i.e.,\begin{eqnarray}
\|S(t)\|_{{\mathcal{L}}{(H^{(1)},H^{(1)})}}\leq C_0^{(1)}.\label{3-8}\end{eqnarray}

\emph{Step 2.} We shall give the proof of the continuity of $S(t)$ with respect to the initial data in $H^{(1)}$.

Assuming  $(u_{0i},v_{0i},\omega_{0i},\theta_{0i})\in H^{(1)}, (u_i,v_i,\omega_i,\theta_i)=S(t)(u_{0i},v_{0i},\omega_{0i},\theta_{0i}),\;(i=1,2)$, and $(u,v,\omega,\theta)=(u_1,v_1,\omega_1,\theta_1)-(u_2,v_2,\omega_2,\theta_2)$.  We subtract the corresponding equations \eqref{1-1} satisfied by $(u_1,v_1,\omega_1,\theta_1)$ and $(u_2,v_2,\omega_2,\theta_2)$, we obtain

\begin{eqnarray}
u_t&=&v_x,\label{3-9}\\
v_t&=&-K\left[\frac{\theta_x}{u_1}-\frac{\theta_1u_x}{u_1^2}+(\frac{1}{u_1}-\frac{1}{u_2})\theta_{2x}
-\left(\frac{\theta_1}{u_1^2}-\frac{\theta_2}{u_2^2}\right)u_{2x}\right]\non\\
&&+\left[\frac{v_{xx}}{u_1}-\frac{v_x}{u_1^2}u_{1x}-\left(\frac{v_{2x}u}{u_1u_2}\right)_x\right],\label{3-10}\\
\omega_t&=&A\left[\frac{\omega_{xx}}{u_1}-\Big(\frac{\omega_{2x}u}{u_1u_2}\Big)_x
-\frac{\omega_xu_{1x}}{u_1^2}-u_1\omega-u\omega_2\right],\label{3-11}\\
C_V\theta_t&=&-K\left(\frac{\theta_1}{u_1}v_x+\left(\frac{\theta_1}{u_1}-\frac{\theta_2}{u_2}\right)v_{2x}\right)
+\frac{v^2_{1x}}{u_1}-\frac{v^2_{2x}}{u_2}+\frac{\omega^2_{1x}}{u_1}-\frac{\omega^2_{2x}}{u_2}\non\\
&&
+u_1(\omega_1^2-\omega_2^2)+u\omega_2^2
+D\left(\frac{\theta_x}{u_1}\right)_x-D\left(\frac{\theta_{2x}}{u_1u_2}\right)_x,\label{3-12}
\end{eqnarray}

and

\begin{eqnarray}
&&t=0:u=u_0,v=v_0,\omega=\omega_0,\theta=\theta_0,\label{3-13}\\
&&x=0,1:v=0,\omega=0,\theta_x=0.\label{3-14}
\end{eqnarray}
By Lemma \ref{3-1}, we know that for any $t>0$ and $i=1,2$, \begin{eqnarray}\label{3-15}
&&\|u_i(t)\|_{H^1}^2+\|v_i(t)\|_{H^1}^2+\|\omega_i(t)\|_{H^1}^2+\|\theta_i(t)\|_{H^1}^2\nonumber\\
&&+\int_0^t(\|u_{ix}\|^2+\|v_i\|_{H^2}^2+\|\omega_i\|_{H^2}^2+\|\theta_{ix}\|_{H^1}^2
+\|v_{it}\|^2+\|\omega_{it}\|^2+\|\theta_{it}\|^2)(s)ds \non\\
&&\leq C_0^{(1)}.\end{eqnarray}

Multiplying \eqref{3-9},  \eqref{3-10},  \eqref{3-11} and  \eqref{3-12} by $u,v,\omega$ and $\theta$, respectively, adding them up and integrating the result over $[0,1]$, and using initial boundary
conditions  \eqref{3-13}-- \eqref{3-14} and  \eqref{3-15}, the embedding theorem and the mean value theorem, we deduce that for any $\varepsilon>0$,

\begin{eqnarray}\label{3-16}
&&\frac{1}{2}\frac{d}{dt}(\|u\|^2+\|v\|^2+\frac{1}{A}\|\omega\|^2+C_V\|\theta\|^2)+\int_0^1\left(\frac{v_x^2}{u_1}
+A\frac{\omega_x^2}{u_1}+D\frac{\theta_x^2}{u_1}\right)dx\nonumber\\
&=&\int_0^1uv_xdx-K\int_0^1\left[\frac{\theta_x}{u_1}-\frac{\theta_1u_x}{u_1^2}+(\frac{1}{u_1}-\frac{1}{u_2})\theta_{2x}
-\left(\frac{\theta_1}{u_1^2}-\frac{\theta_2}{u_2^2}\right)u_{2x}\right]vdx\non\\
&&+\int_0^1\frac{v_{2x}uv_x}{u_1u_2}dx+A\int_0^1\left(\frac{\omega_{2x}u\omega_x}{u_1u_2}-u_1\omega^2-u\omega\omega_2)\right)dx\non\\
&&-K\int_0^1\left(\frac{\theta_1}{u_1}v_x+\left(\frac{\theta_1}{u_1}-\frac{\theta_2}{u_2}\right)v_{2x}\right)\theta dx\non\\
&&+\int_0^1\left(\frac{v^2_{1x}}{u_1}-\frac{v^2_{2x}}{u_2}
+\frac{\omega^2_{1x}}{u_1}-\frac{\omega^2_{2x}}{u_2}
+u_1(\omega_1^2-\omega_2^2)+u\omega_2^2\right)\theta dx\nonumber\\
&&+D\int_0^1\frac{\theta_{2x}\theta_x}{u_xu_2}dx\non\\
&\leq&\varepsilon(\|v_x\|^2+\|\omega_x\|^2+\|\theta_x\|^2)\non\\
&&+C_0^{(1)}
F_1(t)(\|u(t)\|^2+\|v(t)\|_{H^1}^2+\|\omega(t)\|_{H^1}^2+\|\theta(t)\|_{H^1}^2)\end{eqnarray}

where $F_1(t)=1+\|v_{1x}\|_{H^1}^2+\|v_{2x}\|_{H^1}^2+\|\omega_{1x}\|_{H^1}^2+\|\omega_{2x}\|_{H^1}^2+
\|\theta_{1x}\|_{H^1}^2+\|\theta_{2x}\|_{H^1}^2$.

Choosing $\varepsilon$ small enough, by virtue of Lemma \ref{3-1}, we can infer from \eqref{3-16} that

\begin{eqnarray}\label{3-17}
&&\frac{d}{dt}(\|u\|^2+\|v\|^2+\frac{1}{A}\|\omega\|^2+C_V\|\theta\|^2)+1/C_0^{(1)}(\|v_x(t)\|^2+\|\omega_x(t)\|^2+\|\theta_x(t)\|^2)\non\\
&\leq& C_0^{(1)}F_1(t)(\|u(t)\|^2+\|v(t)\|_{H^1}^2+\|\omega(t)\|_{H^1}^2+\|\theta(t)\|_{H^1}^2).
\end{eqnarray}

On the other hand, by Lemma \ref{3-1} and \eqref{3-10}, the embedding theorem and the Cauchy-Schwarz inequality, we get

\begin{eqnarray}
\|v_{xx}(t)\|^2
&\leq&C_0^{(1)}(\|v_t\|^2+\|\theta_x\|^2+\|u_x\|^2+\|v_x\|_{L^{\infty}}^2\|u_{1x}\|^2
+\|v_{2x}\|^2_{L^{\infty}}\|u_x\|^2)\non\\
&\leq&\frac{1}{2}\|v_{xx}(t)\|^2+C_0^{(1)}(\|v_t(t)\|^2+\|v_x\|^2+\|\theta_x\|^2)
+C_0^{(1)}(1+\|v_{2xx}\|^2)\|u\|_{H^1}^2\non\end{eqnarray}
which, leads to
\begin{eqnarray}\label{3-18}
\|v_{xx}(t)\|^2\leq C_0^{(1)}\|v_t(t)\|^2+C_0^{(1)}F_1(t)(\|v_x\|^2+\|u_x\|^2+\|\theta_x\|^2).
\end{eqnarray}

Similarly, we can infer from \eqref{3-11}-\eqref{3-12} that

 \begin{eqnarray}
\|\omega_{xx}(t)\|^2&\leq& C_0^{(1)}\|\omega_t(t)\|^2+C_0^{(1)}F_1(t)(\|u_x\|^2+\|\omega_x\|^2),\label{3-19}\\
\|\theta_{xx}(t)\|^2&\leq& C_0^{(1)}\|\theta_t(t)\|^2+C_0^{(1)}F_1(t)(\|v_x(t)\|^2+\|\theta_x(t)\|^2).\label{3-20}\end{eqnarray}

Differentiating \eqref{3-9} with respect to $x$, multiplying the result by $u_x$ and integrating by parts, and using \eqref{3-18}, we derive that for any small $\epsilon>0$,

\begin{eqnarray}\label{3-21}
\frac{d}{dt}\|u_x(t)\|^2
&=&\int_0^1u_xv_{xx}(x,s)ds\non\\
&\leq&\epsilon\|v_{xx}(t)\|^2+C_0^{(1)}(\epsilon)\|u_x(t)\|^2\non\\
&\leq&C_0^{(1)}\epsilon\|v_t(t)\|^2+C_0^{(1)}(\epsilon)F_1(t)(\|v_x(t)\|^2+\|u_x(t)\|^2+\|\theta_x\|^2).\end{eqnarray}
Multiplying \eqref{3-10} by $v_t$, integrating the results with respect to $x$ on $\Omega$, and using Lemma \ref{le3-1}, \eqref{3-18}, the interpolation inequality and the embedding theorem, we obtain \begin{eqnarray}\label{3-22}
\frac{d}{dt}\int_0^1\frac{v_x^2}{u_1}dx+\frac{1}{C_0^{(1)}}\|v_t(t)\|^2\leq C_0^{(1)}F_1(t)(\|v_x(t)\|^2+\|u_x(t)\|^2+\|\theta_x(t)\|^2).
\end{eqnarray}

Analogously, multiplying \eqref{3-11} and \eqref{3-12} by $\omega_t$ and $\theta_t$, respectively, then integrating the result over $\Omega$, using Lemma \ref{le3-1} and the embedding theorem, we obtain

\begin{eqnarray}
&&\frac{d}{dt}\int_0^1\frac{\omega_x^2}{u_1}dx+\frac{1}{C_0^{(1)}}\|\omega_t(t)\|^2
\leq C_0^{(1)}F_1(t)(\|u_x(t)\|^2+\|\omega_x(t)\|^2),\label{3-23}\\
&&\frac{d}{dt}\int_0^1\frac{\theta_x^2}{u_1}dx+\frac{1}{C_0^{(1)}}\|\theta_t(t)\|^2
\leq C_0^{(1)}F_1(t)(\|v_x(t)\|^2+\|\theta_x(t)\|^2+\|\omega_x\|^2).\label{3-24}\end{eqnarray}
Let \begin{eqnarray}{\mathcal{D}}_1(t)&=&\|u(t)\|^2+\|v(t)\|^2+\|\omega(t)\|^2+\|\theta(t)\|^2
+\|u_x(t)\|^2\nonumber\\
&&+\|\frac{v_x}{\sqrt{u_1}}(t)\|^2+\|\frac{\omega_x}{\sqrt{u_1}}(t)\|^2
+\|\frac{\theta_x}{\sqrt{u_1}}(t)\|^2.\label{3-25}\end{eqnarray}

Adding \eqref{3-17}, \eqref{3-21}, \eqref{3-22}, \eqref{3-23} and \eqref{3-24}, and taking $\epsilon>0$ small enough, we deduce that

\begin{eqnarray}
\frac{d}{dt}{\mathcal{D}}_1(t)&\leq& C_0^{(1)}F_1(t)(\|v_x(t)\|^2+\|\omega_x(t)\|^2+\|u_x(t)\|^2+\|\theta_x(t)\|^2)\non\\
&\leq& C_0^{(1)}F_1(t){\mathcal{D}}_1(t)\label{3-26}\end{eqnarray}
which, by applying  the Gronwall inequality, implies for any $t>0$, \begin{eqnarray}
&&\|u(t)\|_{H^1}^2+\|v(t)\|_{H^1}^2+\|\omega(t)\|_{H^1}^2+\|\theta(t)\|_{H^1}^2\non\\
&\leq& C_0^{(1)}{\mathcal{D}}_{1}(0)\exp\left(C_0^{(1)}\int_0^tF_1(s)ds\right)\non\\
&\leq& C_0^{(1)}\exp(C_0^{(1)}t)(\|u_0\|^2+\|v_0\|^2+\|\omega_0\|^2+\|\theta_0\|^2),\label{3-27}\end{eqnarray}
where, by virtue of \eqref{3-15}, $F_1(t)$ and ${\mathcal{D}}_1(t)$ satisfy \begin{eqnarray}
&&\int_0^tF_1(s)ds\leq C_0^{(1)}\label{3-28}\end{eqnarray} and
\begin{eqnarray}&&\frac{1}{C_0^{(1)}}(\|u(t)\|_{H^1}^2+\|v(t)\|_{H^1}^2+\|\omega(t)\|_{H^1}^2+\|\theta(t)\|_{H^1}^2)
\non\\
&&\leq{\mathcal{D}}_1(t)\leq C_0^{(1)}(\|u(t)\|_{H^1}^2+\|v(t)\|_{H^1}^2+\|\omega(t)\|_{H^1}^2+\|\theta(t)\|_{H^1}^2).\label{3-29}\end{eqnarray}

By \eqref{3-27}, we  know that

\begin{eqnarray}&&\|S(t)(u_{01},v_{01},\omega_{01},\theta_{01})-S(t)(u_{02},v_{02},\omega_{02},\theta_{02})\|_{H^{(1)}}\non\\
&\leq& C_0^{(1)}\exp(C_0^{(1)}t)\|(u_{01},v_{01},\omega_{01},\theta_{01})
-(u_{02},v_{02},\omega_{02},\theta_{02})\|_{H^{(1)}}\label{3-30}\end{eqnarray}
which leads to the continuity of $S(t)$ with respect to  the initial data in $H^{(1)}$.\\

In order to prove \eqref{3-5},  it suffices to show \begin{eqnarray}
S(0)=I\label{3-31}\end{eqnarray}
with $I$ being the unit operator on $H^{(1)}$. To derive \eqref{3-31}, we need show that for any $(u_0,v_0,\omega_0,\theta_0)\in H^{(1)}$, \begin{eqnarray}
\|S(t)(u_0,v_0,\omega_0,\theta_0)-(u_0,v_0,\omega_0,\theta_0)\|_{H^{(1)}}\rightarrow 0, \quad as\quad t\rightarrow 0^+.\label{3-32}\end{eqnarray}

We choose a sequence $(u_0^m,v_0^m,\omega_0^m,\theta_0^m)$
which is smooth enough, for example, $$(u_0^m,v_0^m,\omega_0^m,\theta_0^m)\in (C^{1+\alpha}(\Omega)\times C^{2+\alpha}(\Omega)\times C^{2+\alpha}(\Omega)\times C^{2+\alpha}(\Omega))\cap H^{(1)}$$ for some $\alpha\in (0,1)$, such that \begin{eqnarray}
\|(u_0^m,v_0^m,\omega_0^m,\theta_0^m)-(u_0,v_0,\omega_0,\theta_0)\|_{H^{(1)}}\rightarrow 0,\quad as \quad m\rightarrow +\infty.\label{3-33}\end{eqnarray}

By the regularity results, we can conclude that for arbitrary $T>0$, there exists a unique global smooth solution $$(u^m(t),v^m(t),\omega^m(t),\theta^m(t))\in (C^{1+\alpha}(Q_T)\times C^{2+\alpha}(Q_T)\times C^{2+\alpha}(Q_T)\times C^{2+\alpha}(Q_T))\cap H^{(1)}$$ with $Q_T=\Omega\times (0,T)$.  This gives for  $m=1,2,3,...$
\begin{eqnarray} \|(u^m(t),v^m(t),\omega^m(t),\theta^m(t))-(u_0^m,v_0^m,\omega_0^m,\theta_0^m)\|_{H^{(1)}}\rightarrow 0, \ as\ t\rightarrow 0^+.\label{3-34}\end{eqnarray}

Fixing $T=1$, by the continuity of the operator $S(t)$, \eqref{3-32} and \eqref{3-34}, for any $t\in[0,1]$,\begin{eqnarray}
&&\hspace{-0.8cm}\|(u^m(t),v^m(t),\omega^m(t),\theta^m(t))-(u(t),v(t),\omega(t),\theta(t))\|_{H^{(1)}}\non\\
&=&\|S(t)(u_0^m,v_0^m,\omega_0^m,\theta_0^m)-S(t)(u_0,v_0,\omega_0,\theta_0)\|_{H^{(1)}}\non\\
&\leq& C_0^{(1)}\|(u_0^m,v_0^m,\omega_0^m,\theta_0^m)-(u_0,v_0,\omega_0,\theta_0)\|_{H^{(1)}}
\rightarrow 0,\quad as \quad m\rightarrow +\infty.\non\end{eqnarray}

Thus, this along with  \eqref{3-33} and \eqref{3-34}, gives
\begin{eqnarray}
&&\hspace{-0.8cm}\|S(t)(u_0,v_0,\omega_0,\theta_0)-(u_0,v_0,\omega_0,\theta_0)\|_{H^{(1)}}\non\\
&=&\|(u(t),v(t),\omega(t),\theta(t)-(u_0,v_0,\omega_0,\theta_0)\|_{H^{(1)}}\non\\
&\leq&\|(u^m(t),v^m(t),\omega^m(t),\theta^m(t))-(u(t),v(t),\omega(t),\theta(t))\|_{H^{(1)}}\non\\
&&
+\|(u^m(t),v^m(t),\omega^m(t),\theta^m(t))-(u_0^m,v_0^m,\omega_0^m,\theta_0^m)\|_{H^{(1)}}\non\\
&&
+\|(u_0^m,v_0^m,\omega_0^m,\theta_0^m)-(u_0,v_0,\omega_0,\theta_0)\|_{H^{(1)}}\rightarrow 0,\; as \; m\rightarrow +\infty,\;t\rightarrow 0^+,\non
\end{eqnarray}
which implies \eqref{3-31} and  \eqref{3-32}.  By  \eqref{3-7},  \eqref{3-8} and  \eqref{3-31}. Hence, we conclude $S(t)$ is a $C_0$-semigroup on $H^{(1)}$ satisfying \eqref{3-5}.\\

For any $t_1>0$, integrating the third equation of \eqref{1-1} over $(t_1,t)$ and using Lemma \ref{3-1}, we obtain \begin{eqnarray}
\|\omega(t)-\omega(t_1)\|&\leq& C_0^{(1)}\left|\int_{t_1}^t(\|\omega_{xx}\|^2+\|w_x\|_{L^{\infty}}^2\|u_x\|^2+\|\omega\|^2)ds\right|^{\frac{1}{2}}
\left|t-t_1\right|^{\frac{1}{2}}\non\\
&\leq&C_0^{(1)}\left|\int_{t_1}^t(\|\omega_{xx}\|^2+\|w_x\|^2+\|u_x\|^2)ds\right|^{\frac{1}{2}}
\left|t-t_1\right|^{\frac{1}{2}}\non\\
&\leq& C_0^{(1)}\left|t-t_1\right|^{\frac{1}{2}}\non\end{eqnarray}
which implies $$ \omega(t)\in C^{1/2}([0,+\infty),L^2).$$
In the same manner, we can prove $u(t)\in C^{1/2}([0,+\infty),H^1), v(t),\theta(t)\in C^{1/2}([0,+\infty),L^2).$ Thus, we can obtain \eqref{3-6}.
The proof is complete.
\hspace*{\fill}$\Box$

\begin{Lemma}\label{le3-4} Under the assumptions in Theorem \ref{th2-1}, the problem \eqref{1-1}-\eqref{1-3} admits a unique generalized global solution $(u(t),v(t),\omega(t),\theta(t))$ in $H^{(2)}$ which
	defines a nonlinear $C_0$-semigroup $S(t)$ (also denoted by $S(t)$ by the uniqueness of solution in $H^{(1)}$) on $H^{(2)}$ such that for any
	$(u_0,v_0,\omega_0,\theta_0)\in H^{(2)}$, the generalized global solution $(u(t),v(t),\omega(t),\theta(t))$ satisfies
	\begin{eqnarray}
	&&\|S(t)(u_0,v_0,\omega_0,\theta_0)\|_{H^{(2)}}=\|(u(t),v(t),\omega(t),\theta(t))_{H^{(2)}}\leq C_0^{(2)},\label{3-35}\\
	&&(u(t),v(t),\omega(t),\theta(t))=S(t)(u_0,v_0,\omega_0,\theta_0)\in C([0,\infty), H^{(2)}),\label{3-36}\\
	&&u(t)\in C^{\frac{1}{2}}([0,+\infty),H^2),\quad v(t),\omega(t),\theta(t)C^{\frac{1}{2}}([0,+\infty),H_1).\label{3-37}
	\end{eqnarray}
\end{Lemma}

\noindent {\bf Proof. } The estimate \eqref{3-35} and the global existence of generalized solution $(u(t),v(t),\omega(t),\\\theta(t))\in H^{(2)}$
follow from Lemma \ref{le3-2}. Similarly to Lemma \ref{3-3}, we can prove the estimate \eqref{3-37}.
In order to complete the proof of Lemma \ref{le3-4}, it suffices to prove \eqref{3-36} and the continuity
of $S(t)$  with respect to $(u_0,v_0,\omega_0,\theta_0)$ in $H^{(2)}$, which also leads to the uniqueness of the generalized global solutions in $H^{(2)}$. This will be done as follows.\\

The uniqueness of generalized global solutions in $H^{(2)}$ follows from that in $H^{(1)}$. Thus
$S(t)$ satisfies \eqref{3-7} on $H^{(2)}$ and by Lemma \ref{le3-2}, \begin{eqnarray}
\|S(t)\|_{{\mathcal{L}}{(H^{(2)},H^{(2)})}}\leq C_0^{(2)}.\label{3-38}\end{eqnarray}
In the same manner as in the proof of Lemma \ref{le3-3}, we
assume that  $(u_{0i},v_{0i},\omega_{0i},\theta_{0i})\in H^{(2)}$, $(u_i,v_i,\omega_i,\theta_i)=S(t)(u_{0i},v_{0i},\omega_{0i},\theta_{0i}),\;(i=1,2)$, and $(u,v,\omega,\theta)=(u_1,v_1,\omega_1,\theta_1)-(u_2,v_2,\omega_2,\theta_2)$.  We subtract the corresponding equations \eqref{1-1} satisfied by $(u_1,v_1,\omega_1,\theta_1)$ and $(u_2,v_2,\omega_2,\theta_2)$, we obtain equations \eqref{3-9}-\eqref{3-12}.\\

By Lemma \ref{3-2}, we know that for any $t>0$ and $i=1,2$, \begin{eqnarray}
&&\|u_i(t)\|_{H^2}^2+\|v_i(t)\|_{H^2}^2+\|\omega_i(t)\|_{H^2}^2+\|\theta_i(t)\|_{H^2}^2
+\|v_{it}(t)\|^2+\|\omega_{it}(t)\|^2+\|\theta_{it}(t)\|^2\non\\
&&\qquad
+\int_0^t(\|u_{ix}\|_{H^1}^2+\|v_i\|_{H^3}^2+\|\omega_i\|_{H^3}^2+\|\theta_{ix}\|_{H^2}^2
+\|v_{it}\|_{H^1}^2+\|\omega_{it}\|_{H^1}^2+\|\theta_{it}\|_{H^1}^2)(s)ds\non\\
&&\leq C_0^{(2)}.\label{3-39}\end{eqnarray}

By \eqref{3-18}-\eqref{3-20}, we have

\begin{eqnarray} \|v_{xx}(t)\|^2&\leq& C_0^{(2)}(\|v_t(t)\|^2+\|v_x(t)\|^2+\|u_x(t)\|^2+\|\theta_x(t)\|^2),\label{3-40}\\
\|\omega_{xx}(t)\|^2&\leq& C_0^{(2)}(\|\omega_t(t)\|^2+\|u_x(t)\|^2+\|\omega_x(t)\|^2),\label{3-41}\\
\|\theta_{xx}(t)\|^2&\leq& C_0^{(2)}(\|\theta_t(t)\|^2+\|v_x(t)\|^2+\|\theta_x(t)\|^2),\label{3-42}.\end{eqnarray}
By virtue of  \eqref{3-39},  $$F_1(t)=1+\|v_{1x}\|_{H^1}^2+\|v_{2x}\|_{H^1}^2+\|\omega_{1x}\|_{H^1}^2+\|\omega_{2x}\|_{H^1}^2+
\|\theta_{1x}\|_{H^1}^2+\|\theta_{2x}\|_{H^1}^2\leq C_0^{(2)}.$$
Differentiating \eqref{3-10} with respect to $x$, we have \begin{eqnarray}
v_{tx}=\frac{v_{xxx}}{u_1}-\frac{2v_{xx}u_{1x}}{u_1^2}+M(x,t),\label{3-43}\end{eqnarray}
where \begin{eqnarray} M(x,t)
&=&-K\bigg[\frac{\theta_{xx}}{u_1}-\frac{\theta_xu_{1x}}{u_1^2}-\frac{\theta_{1x}u_x
	+\theta_1u_{xx}}{u_1^2}
+\frac{2\theta_1u_xu_{1x}}{u_1^3}+(\frac{1}{u_1}-\frac{1}{u_2})\theta_{2xx}\non\\
&&
+\left(\frac{u_{2x}}{u_2^2}-\frac{u_{1x}}{u_1^2}\right)\theta_{2x}
-\left(\frac{\theta_1}{u_1^2}-\frac{\theta_2}{u_2^2}\right)u_{2xx}-\Big(\frac{\theta_{1x}}{u_1^2}
-\frac{\theta_{2x}}{u_2^2}\non\\
&&-\frac{2\theta_1u_{1x}}{u_1^3}
+\frac{2\theta_2u_{2x}}{u_2^3}\Big)u_{2x}\bigg]-\frac{v_xu_{1xx}}{u_1^2}+\frac{2v_xu_{1x}^2}{u_1^3}
-\left(\frac{v_{2x}u}{u_1u_2}\right)_{xx}.\non\end{eqnarray}

By Lemmas \ref{le3-1} and \ref{le3-2}, \eqref{3-39}, the embedding theorem and the Gagliardo-Nirenberg interpolation inequality, we obtain

\begin{eqnarray}
\|M(t)\|^2
&\leq& C_0^{(2)}(\|\theta_{xx}\|^2+\|\theta_x\|^2\|u_{1x}\|_{L^{\infty}}^2
+\|u_x\|^2\|\theta_{1x}\|_{L^{\infty}}^2+\|u_{xx}\|^2+\|u_x\|^2\|u_{1x}\|_{L^{\infty}}^2
\non\\
&&+\|\theta_{2xx}\|^2+\|\theta_{2x}\|^2\|u_{1x}\|_{L^{\infty}}^2
+\|\theta_{2x}\|^2\|u_{2x}\|_{L^{\infty}}^2
+\|u_{2xx}\|^2+\|\theta_{1x}\|^2\|u_{2x}\|_{L^{\infty}}^2\non\\
&&
+\|\theta_{2x}\|^2\|u_{2x}\|_{L^{\infty}}^2
+\|u_{2x}\|^2\|u_{2x}\|_{L^{\infty}}^2+\|u_{1x}\|^2\|u_{2x}\|_{L^{\infty}}^2
+\|u_{1xx}\|^2\|v_{x}\|_{L^{\infty}}^2\non\\
&&+\|v_{x}\|^2\|u_{1x}\|_{L^{\infty}}^4
+\|v_{2xxx}\|^2\|u\|_{L^{\infty}}^2+\|u_{x}\|^2\|v_{2xx}\|_{L^{\infty}}^2
+\|u_{xx}\|^2\|v_{2x}\|_{L^{\infty}}^2
\non\\
&&+\|v_{2xx}\|^2\|u_{1x}\|_{L^{\infty}}^2+\|v_{2xx}\|^2\|u_{2x}\|_{L^{\infty}}^2
+\|v_{2x}\|_{L^{\infty}}^2\|u_{1x}\|_{L^{\infty}}^2\|u_x\|^2\non\\
&&
+\|v_{2x}\|_{L^{\infty}}^2\|u_{2x}\|_{L^{\infty}}^2\|u_x\|^2)\non\\
&\leq& C_0^{(2)}(1+\|v_{2xxx}\|^2)(\|v_x\|_{H^1}^2+\|\theta_x\|^2_{H^1}+\|u_x\|^2).
\label{3-44}\end{eqnarray}

By \eqref{3-43}, \eqref{3-44} and the interpolation inequality, we can infer that \begin{eqnarray}
\|v_{xxx}\|^2&\leq& C_0^{(1)}\|v_{tx}\|^2+C_0^{(2)}(\|v_{xx}\|_{L^{\infty}}^2+\|M(t)\|^2)\non\\
&\leq&\frac{1}{2}\|v_{xxx}\|^2+C_0^{(1)}\|v_{tx}\|^2\nonumber\\
&&+C_0^{(2)}(1+\|v_{2xxx}\|^2)(\|v_x\|_{H^1}^2+\|\theta_x\|^2_{H^1}+\|u_x\|^2)\non\end{eqnarray}
which, gives \begin{eqnarray} \|v_{xxx}\|^2\leq C_0^{(1)}\|v_{tx}\|^2+C_0^{(2)}(1+\|v_{2xxx}\|^2)(\|v_x\|_{H^1}^2+\|\theta_x\|^2_{H^1}+\|u_x\|^2).\label{3-45}\end{eqnarray}

Differentiating \eqref{3-9} twice with respect to $x$, multiplying the result by $u_{xx}$, integrating the result equation over $\Omega$, using \eqref{3-45} and the Cauchy inequality, we have \begin{eqnarray}
\frac{d}{dt}\|u_{xx}(t)\|^2&\leq& C_0^{(1)}(\|u_{xx}(t)\|^2+\|v_{xxx}(t)\|^2)\non\\
&\leq& C_0^{(1)}\|v_{tx}\|^2+C_0^{(2)}(1+\|v_{2xxx}\|^2)(\|v_x\|_{H^1}^2+\|\theta_x\|^2_{H^1}+\|u_x\|_{H^1}^2).\label{3-46}\end{eqnarray}

Differentiating \eqref{3-10} with respect to $t$, multiplying the result by $v_t$, integrating the result over $\Omega$, using Lemmas \ref{le3-1}-\ref{le3-3} and \eqref{3-39}, we deduce that \begin{eqnarray}
\frac{d}{dt}\|v_t\|^2+\frac{1}{C_0^{(1)}}\|v_{tx}\|^2&\leq& C_0^{(2)}(1+\|v_{2tx}\|^2+\|\theta_{1tx}\|^2)(\|v_x\|_{H^1}^2+\|\theta_x\|^2\non\\
&&+\|u_x\|^2
+\|v_t\|^2+\|\theta_t\|^2).\label{3-47}\end{eqnarray}

Similarly to \eqref{3-47}, differentiating \eqref{3-11} and \eqref{3-12} with respect to $t$, multiplying it by $\omega_t$ and $\theta_t$, respectively,  using the embedding theorem and Lemmas \ref{le3-1}-\ref{le3-3} and \eqref{3-39}, we have for any small $\varepsilon>0$, we have \begin{eqnarray}\frac{d}{dt}\|\omega_t\|^2+\frac{1}{C_0^{(1)}}\|\omega_{tx}\|^2\leq C_0^{(2)}(1+\|\omega_{2tx}\|^2)(\|v_x\|^2+\|\omega_t\|^2+\|\omega\|^2_{H^1}+\|u\|^2).\label{3-48}\end{eqnarray}
and \begin{eqnarray}
&&\frac{C_V}{2}\frac{d}{dt}\|\theta_t\|^2+D\int_0^1\frac{\theta_{tx}^2}{u_x}dx\non\\
&=&D\int_0^1\left[\left(\frac{\theta_{2x}}{u_1u_2}\right)_t
+\frac{\theta_xv_{1x}}{u_1^2}\right]\theta_{tx}dx+\int_0^1\left(\frac{v_{1x}^2}{u_1}
-\frac{v_{2x}^2}{u_2}\right)_t\theta_tdx\non\\
&&-K\int_0^1\left[\frac{\theta_1v_x}{u_1}
+\left(\frac{\theta_1}{u_1}-\frac{\theta_2}{u_2}\right)v_{2x}\right]_t\theta_tdx
+\int_0^1\left(\frac{\omega_{1x}^2}{u_1}-\frac{\omega_{2x}^2}{u_2}\right)_t\theta_tdx\non\\
&&
+\int_0^1\left(u_1(\omega_1^2-\omega_2^2)+u\omega_2^2\right)_t\theta_tdx\non\\
&\leq& \varepsilon\|\theta_{tx}\|^2+ C_0^{(2)}\|v_{tx}\|^2+C_0^{(2)}(1+\|\theta_{1tx}\|^2+\|v_{2tx}\|^2+\|\omega_{1tx}\|^2+\|\omega_{2tx}\|^2)\non\\
&&\times(\|\theta_t\|^2+\|u_x\|^2+\|v_x\|^2+\|v_t\|^2+\|\theta_x\|^2).\label{3-49}\end{eqnarray}
Choosing $\varepsilon$ small enough, we derive from \eqref{3-49} that \begin{eqnarray}
\frac{d}{dt}\|\theta_t\|^2+C_0^{(1)}\|\theta_{tx}\|^2
&\leq& C_0^{(2)}\|v_{tx}\|^2+C_0^{(2)}(1+\|\theta_{1tx}\|^2+\|v_{2tx}\|^2+\|\omega_{1tx}\|^2+\|\omega_{2tx}\|^2)\non\\
&&\times(\|\theta_t\|^2+\|u_x\|^2+\|v_x\|^2+\|v_t\|^2+\|\theta_x\|^2).\label{3-50}\end{eqnarray}
Now multiplying \eqref{3-47} by a large number $N>2(C_0^{(1)})^2$, then adding up the result, \eqref{3-46}, \eqref{3-48} and \eqref{3-50}, we conclude\begin{eqnarray}
\frac{d}{dt}{\mathcal{D}}_2(t)&\leq& C_0^{(2)}F_2(t)(\|v_t\|^2+\|\omega_t\|^2+\|\theta_x\|^2+\|v_x\|_{H^1}^2+\|\theta_x\|^2_{H^1}+\|u_x\|_{H^1}^2)\non\\
&\leq& C_0^{(2)}F_2(t)({\mathcal{D}}_1(t)+{\mathcal{D}}_2(t))\label{3-51}\end{eqnarray}
where ${\mathcal{D}}_2(t)=\|u_{xx}\|^2+N\|v_t\|^2+\|\omega_t\|^2+\|\theta_t\|^2$ and $F_2(t)=1+\|v_{2xxx}\|^2+\|\theta_{1tx}\|^2+\|\theta_{2tx}\|^2+\|v_{1tx}\|^2+\|v_{2tx}\|^2
+\|\omega_{1tx}\|^2+\|\omega_{2tx}\|^2$.\\
By \eqref{3-39}, $F_2(t)$ satisfies \begin{eqnarray}\int_0^tF_2(s)ds\leq C_0^{(2)}(1+t),\quad \forall t>0.\label{3-52}\end{eqnarray}
Adding \eqref{3-26} to \eqref{3-51} gives \begin{eqnarray}
\frac{d}{dt}{\mathcal{D}}(t)\leq C_0^{(2)}F_2(t){\mathcal{D}}(t),\label{3-53}\end{eqnarray}
where, by \eqref{3-18}-\eqref{3-20}, \eqref{3-40}-\eqref{3-42}, ${\mathcal{D}}(t)={\mathcal{D}}_1(t)+{\mathcal{D}}_2(t)$ satisfies\begin{eqnarray}
&&\frac{1}{C_0^{(2)}}(\|u(t)\|_{H^2}^2+\|v(t)\|_{H^2}^2+\|\omega(t)\|_{H^2}^2+\|\theta(t)\|_{H^2}^2)
\non\\
&&\leq {\mathcal{D}}(t)\leq C_0^{(2)}(\|u(t)\|_{H^2}^2+\|v(t)\|_{H^2}^2+\|\omega(t)\|_{H^2}^2+\|\theta(t)\|_{H^2}^2).\label{3-54}\end{eqnarray}
Thus it follows from \eqref{3-53}, Gronwall's inequality, and \eqref{3-54} that \begin{eqnarray}
&&\|u(t)\|_{H^2}^2+\|v(t)\|_{H^2}^2+\|\omega(t)\|_{H^2}^2+\|\theta(t)\|_{H^2}^2\non\\
&\leq& C_0^{(2)}{\mathcal{D}}(t)\leq C_0^{(2)}{\mathcal{D}}(0)\exp\left(C_0^{(2)}\int_0^tF_2(s)ds\right)\non\\
&\leq&\exp(C_0^{(2)}t)(\|u_0\|^2_{H^1}+\|v_0\|^2_{H^1}+\|\omega_0\|^2_{H^1}
+\|\theta_0\|^2_{H^1}),\ \forall t>0,\label{3-55}\end{eqnarray}
which, implies the continuity of $S(t)$ with respect to the initial data in $H^{(2)}$. Similarly to
the proof of \eqref{3-5}, we can prove \eqref{3-36}. The proof is done. \hspace*{\fill}$\Box$

\section{Attracting property}
\setcounter{equation}{0}
\subsection{ Existence of an absorbing set in $H_{\delta}^{(1)}$}\label{sec4}
In this section, we shall show the existence of an absorbing set in $H_{\delta}^{(1)}$. Throughout this section we always assume that the initial data belong to a bounded set of $H_{\delta}^{(1)}$. We begin with the following lemma.
\begin{Lemma}\label{le4-1}
	If $(u_0,v_0,\omega_0,\theta_0)\in H_{\delta}^{(1)}$, then the following estimates hold
	for any $(x,t)\in [0,1]\times [0,+\infty)$:
	\begin{eqnarray}\label{4-1}\begin{cases}
	\delta_3\leq\int_0^1u(x,t)dx=\int_0^1u_0(x)dx\leq\delta_4,\;\;\forall t>0,&\\
	\delta_5\leq\int_0^1(\frac{1}{2}v^2+\frac{1}{2A}\omega^2+C_V\theta)dx\leq\delta_2,&\\
	-\int_0^1(K\log u+C_V\log\theta)(x,t)dx+\int_0^t\int_0^1\left(\frac{v_x^2}{u\theta}+\frac{\omega_x^2}{u\theta}
	+\frac{u\omega^2}{\theta}+D\frac{\theta_x^2}{u\theta^2}\right)dxds\leq-\delta_1,&\\
	0<C_{\delta}^{-1}\leq u(x,t)\leq C_{\delta},&\\
	\theta(x,t)\geq C_{\delta}^{-1}>0.&\end{cases}\end{eqnarray}
\end{Lemma}
{\bf Proof. }
See, e.g., \cite{HN} and \cite{M3}.\hspace*{\fill}$\Box$

\begin{Lemma}\label{le4-2}
	For initial data belonging to a bounded set of $H_{\delta}^{(1)}$, there is $t_0>0$,
	depending only on the boundedness of this set, such that for all $t\geq t_0,\;x\in[0,1]$,
	\begin{eqnarray}\label{4-2}\frac{\delta_3}{2}\leq u(x,t)\leq 2\delta_4,\quad \frac{\delta_5}{2C_V}\leq \theta(x,t)\leq \frac{2\delta_2}{C_V}.
	\end{eqnarray}
\end{Lemma}
{\bf Proof. }  It can be seen from the first equation of \eqref{1-1} and \eqref{3-2} that as $t\rightarrow+\infty$,
\begin{eqnarray}\int_0^1u(x,t)dx=\int_0^1u_0(x)dx,\ \ \left\|u-\int_0^1udx\right\|_{L^\infty}\leq C_0^{(1)}\left\|u-\int_0^1udx\right\|_{H^1}\rightarrow 0.\label{4-3}\end{eqnarray}
we will use a contradiction argument to prove \eqref{4-2}. Suppose that it is not true, then there exists a sequence $\{t_n\}\uparrow+\infty$ such that for all $x\in [0,1]$,\begin{eqnarray}\sup u(x,t_n)>2\delta_4,\label{4-4}\end{eqnarray}
where the sup is taken for all initial data in a given bounded set of $H_{\delta}^{(1)}$. Then in the same manner as the proof of Lemma 3.3 in  \cite{zq1}, there is $(u_0,v_0,\omega_0,\theta_0)$ belonging to this bounded set such that for the corresponding solution $(u,v,\omega,\theta)$, we have \begin{eqnarray}
u(x,t_n)\geq 2\delta_4,\quad \forall x\in[0,1].\label{4-5}\end{eqnarray}
This contradicts with \eqref{4-3} and \eqref{4-1}. In the same way, we can derive other parts of \eqref{4-2}.  \hspace*{\fill}$\Box$\\

It follows from Lemmas \ref{le4-1} and \ref{le4-2} that, for initial data belonging to a given bounded set ${\mathcal{B}}$ of $H^{(1)}_{\delta}$, the orbit will re-enter $H^{(1)}_{\delta}$ and stay there after a finite time. In the following, we shall prove the existence of an absorbing ball in $H^{(1)}_{\delta}$. Since we assume that the initial data $(u_0,v_0,\omega_0,\theta_0)\in {\mathcal{B}}_1$( ${\mathcal{B}}_1$ is an arbitrary bounded set of $H^{(1)}_{\delta}$), there is a positive constant $B_1$ such that $\|(u_0,v_0,\omega_0,\theta_0)\|_{H^1}\leq B_1$, and we use $C_{B_1,\delta}$ to denote generic positive constant depending on $B_1$ and $\delta_i(i=1,2,3,4,5)$.

\begin{Lemma}\label{le4-3}
	There exists a positive constant $\hat{\gamma}_1=\hat{\gamma}_1(C_{B_1,\delta})>0$
	such that, for any fixed $\gamma\in(0,\hat{\gamma}_1]$, it holds that for any $t>0$,\begin{eqnarray}
	&& e^{\gamma t}(\|u-\bar
	u\|^2_{H^1}+\|v\|_{H^1}^2+\|\omega\|_{H^1}^2+\|\theta-\bar
	\theta\|^2_{H^1})\non\\
	&&+\int_0^te^{\gamma s}(\|u-\bar
	u\|^2_{H^1}+\|v\|_{H^2}^2+\|\theta-\bar
	\theta\|^2_{H^2}+\|\omega\|_{H^2}^2+\|v_t\|^2+\|\omega_t\|^2+\|\theta_t\|^2)(s)\,d\tau\non\\
	&&\leq
	C_{B_1,\delta}\label{4-6}\end{eqnarray} which implies
	\begin{eqnarray}
	\|u(t)\|_{H^1}^2+\|v(t)\|_{H^1}^2+\|\omega(t)\|_{H^1}^2+\|\theta(t)\|_{H^1}^2&\leq& 2({\bar u}^2+{\bar\theta}^2)+C_{B_1,\delta}e^{-\gamma t}\non\\
	& \leq& 2\left(\delta_4^2+\frac{\delta_2^2}{C_V^2}\right)+C_{B_1,\delta}e^{-\gamma t}.\label{4-7}
	\end{eqnarray}
\end{Lemma}
{\bf Proof }   By virtue of Theorem 1.1 in \cite{HN}, we can claim the same argument and  easily prove this lemma. \hspace*{\fill}$\Box$

Therefore, the following results on the existence of an absorbing  set in $H^{(1)}_{\delta}$ follow from Lemma \ref{le4-3}.

\begin{Lemma}\label{le4-4}
	Let\begin{eqnarray} R_1(\delta)=2\sqrt{\frac{C_V^2\delta_4^2+\delta_2^2}{C_V^2}},\quad and \quad {\mathcal{B}}_1=\{(u,v,\omega,\theta)\in H_{\delta}^{(1)},\|(u,v,\omega,\theta)\|_{H^{(1)}}\leq R_1\}\non.\end{eqnarray} Then ${\mathcal{B}}_1$ is an absorbing ball in  $H_{\delta}^{(1)}$, i.e., there exists some $$t_1=t_1(C_{B,\delta})=\max\left\{-\hat\gamma_1^{-1}
	\log[2\frac{C_V^2\delta_4^2+\delta_2^2}{C_{B,\delta}C_V^2}],t_0\right\}\geq t_0$$ such that  when $t\geq t_1$, $$\|(u(t),v(t),\omega(t),\theta(t))\|_{H^{(1)}}^2\leq R_1^2.$$
\end{Lemma}

\subsection{ Existence of an absorbing set in $H_{\delta}^{(2)}$}\label{sec5}
\setcounter{equation}{0} In this section, we address the existence of an absorbing set in $H_{\delta}^{(2)}$. Throughout this section we assume that the initial data is in an arbitrarily fixed bounded set ${\mathcal{B}}_2$ in $H_{\delta}^{(2)}$, i.e., $\|(u_0,v_0,\omega_0,\theta_0)\|_{H^2}\leq B_2$ with $B_2,$ a given a positive constant.
\begin{Lemma}\label{le5-1}
	There exists a positive constant $\hat{\gamma}_2=\hat{\gamma}_2(C_{B_2,\delta})\leq\hat{\gamma}_1$
	such that, for any fixed $\gamma\in(0,\hat{\gamma}_2]$, it holds that for any $t>0$,\begin{eqnarray}
	&& e^{\gamma t}(\|u-\bar
	u\|^2_{H^2}+\|v\|_{H^2}^2+\|\omega\|_{H^2}^2+\|\theta-\bar
	\theta\|^2_{H^2})\non\\
	&&+\int_0^te^{\gamma s}(\|u-\bar
	u\|^2_{H^2}+\|v\|_{H^3}^2+\|\theta-\bar
	\theta\|^2_{H^3}+\|\omega\|_{H^3}^2+\|v_t\|_{H^1}^2+\|\omega_t\|_{H^1}^2+\|\theta_t\|_{H^1}^2)(s)\,ds\nonumber\\
&&\leq
	C_{B_2,\delta}\label{5-1}\end{eqnarray} which implies
	\begin{eqnarray}
	\|u(t)\|_{H^1}^2+\|v(t)\|_{H^1}^2+\|\omega(t)\|_{H^1}^2+\|\theta(t)\|_{H^1}^2
	&\leq& 2({\bar u}^2+{\bar\theta}^2)+C_{B_2,\delta}e^{-\gamma t}\non\\
	&\leq& 2\left(\delta_4^2+\frac{\delta_2^2}{C_V^2}\right)+C_{B_2,\delta}e^{-\gamma t}.\label{5-2}
	\end{eqnarray}
\end{Lemma}
{\bf Proof. }  The proof is similar as in \cite{HN}, and we omit the detail for this lemma.  \hspace*{\fill}$\Box$\\

By Lemma \ref{le5-1}, we immediately obtain the following Lemma:

\begin{Lemma}\label{le5-2}
	Let\begin{eqnarray} R_2(\delta)=2\sqrt{\frac{C_V^2\delta_4^2+\delta_2^2}{C_V^2}}\non\end{eqnarray}
	and \begin{eqnarray} {\mathcal{B}}_2=\{(u,v,\omega,\theta)\in H_{\delta}^{(2)},\|(u,v,\omega,\theta)\|_{H^{(2)}}\leq R_1\}\non.\end{eqnarray} Then ${\mathcal{B}}_2$ is an absorbing ball in  $H_{\delta}^{(2)}$, i.e., there exists some $$t_2=t_2(C_{B_2,\delta})\geq\max\left\{-\hat\gamma_2^{-1}
	\log\left[2\frac{C_V^2\delta_4^2+\delta_2^2}{C_{B_2,\delta}C_V^2}\right],t_1(C_{B_1,\delta})\right\}$$ such that  when $t\geq t_2$, $$\|(u(t),v(t),\omega(t),\theta(t))\|_{H^{(2)}}^2\leq R_2^2.$$
\end{Lemma}

\section{Proof of main result}\label{sec6}
\setcounter{equation}{0}
\subsection{Preliminary theory of global attractor}
\begin{Lemma}\label{le6-1}
	Let $H_1,\;H_2,\;H_3$ be three Banach spaces verifying the following conditions:\\
	
	\begin{itemize}
	(1) the embedding $H_3\hookrightarrow H_2$ and $H_2 \hookrightarrow H_1$ are compact;\\
	(2) there exists a $C_0$-semigroup $\{S(t)\}$ on $H_2$ and $H_3$ which maps $H_2, H_3$ into $H_2$ and $H_3$, respectively,
	and for any $t>0$, $S(t)$ is continuous (nonlinear) operator on $H_2$ and $H_3$, respectively;\\
	(3) the semigroup $S(t)$ on $H_3$ possesses a bounded absorbing set in $H_3$; \\
	
	then there is a weak universal attractor ${\mathcal{A}}_3$ in $H_3$.\\
	
	If, furthermore, the following conditions are satisfied:\\
	
	(4) the semigroup $S(t)$ on $H_2$ possesses a bounded absorbing set in $H_2$;\\
	(5) for any $t>0$, $S(t)$ is continuous on bounded sets of $H_2$ for the topology the norm of $H_1$;\\
	
	then there is a weak universal attractor ${\mathcal{A}}_2$ in $H_2$.
        \end{itemize}
\end{Lemma}
{\bf Proof.}  See, e.g., Ghidaglia \cite{6}. \hspace*{\fill}$\Box$

\subsection{Some lemmas to construct $\omega$-limit set}
We have proved the existence of absorbing balls in $H^{(1)}_{\delta}$ in Section \ref{sec4}, then we can use Lemma \ref{le6-1} to prove Theorem \ref{th2-1}.

\begin{Lemma}\label{le6-2}
	The set \begin{eqnarray}
	\omega({\mathcal{B}}_2)=\bigcap_{s\geq0}\overline{\bigcup_{t\geq s}S(t){\mathcal{B}}_2},\label{6-1}\end{eqnarray}
	where the closures are taken with respect to the weak topology of $H^{(2)}_{\delta}$, is included in ${\mathcal{B}}_2$ and is nonempty. It is invariant under operators $\{S(t)\}$, i.e.,\begin{eqnarray}
	S(t)\omega({\mathcal{B}}_2)=\omega({\mathcal{B}}_2),\quad \forall t>0.\label{6-2}\end{eqnarray}
\end{Lemma}

\begin{Lemma}\label{le6-3}
	The set \begin{eqnarray}
	{\mathcal{A}}_{2,\delta}=\omega({\mathcal{B}}_2)\label{6-3}\end{eqnarray}
	satisfies\begin{eqnarray} {\mathcal{A}}_{2,\delta}\quad  is\;\; bounded\;\; and\;\; weakly\;\; closed\;\; in H^{(2)}_{\delta},\
	S(t){\mathcal{A}}_{2,\delta}={\mathcal{A}}_{2,\delta},\quad \forall t\geq 0,\label{6-4}\end{eqnarray}
	and, for every bounded set ${\mathcal{B}}$ in $H_{\delta}^{(2)}$,\begin{eqnarray}
	\lim_{t\rightarrow+\infty}d^{\omega}(S(t){\mathcal{B}},{\mathcal{A}}_{2,\delta})=0.\label{6-5}\end{eqnarray} Moreover, it is the maximal set in the sense of inclusion that satisfies \eqref{6-4}-\eqref{6-6}.
\end{Lemma}
{\bf Proof of Lemma \ref{le6-2} and \ref{le6-3}. }  This proof follows from Lemma \ref{le6-1}, due to the fact that $S(t)$ is continuous on $H^{(2)}_{\delta}$ and $H^{(1)}_{\delta}$, respectively, and $H^{(2)}_{\delta}$ is compactly embedded in $H^{(1)}_{\delta}$, ${\mathcal{B}}_2$ and ${\mathcal{B}}_1$ are absorbing balls in $H^{(2)}_{\delta}$ and $H^{(2)}_{\delta}$, respectively. \hspace*{\fill}$\Box$\\

Leading similar fashion in \cite{6}, we also call ${\mathcal{A}}_{2,\delta}$ the universal attractor of $S(t)$ in $H_{\delta}^{(2)}$. In order to discuss the existence of a universal attractor in $H_{\delta}^{(1)}$, we need to prove the following lemma:

\begin{Lemma}\label{le6-4}
	For every $t\geq0$, operator $S(t)$ is continuous in $H_{\delta}^{(1)}$ for the topology induced by the norm of $L^2\times L^2\times L^2\times L^2$.
\end{Lemma}
{\bf Proof.} The proof is the same as in Lemma \ref{le3-3}, we can repeat the same argument as the proof of \eqref{3-26} in $H^{(1)}_{\delta}$, and complete the proof of this lemma.  \hspace*{\fill}$\Box$\\

Now we can again use Lemma \ref{le6-1} to obtain the following result on existence of a universal attractor in $H^{(1)}_{\delta}$.

\begin{Lemma}\label{le6-5}
	The set \begin{eqnarray}
	{\mathcal{A}}_{1,\delta}=\bigcap_{s\geq0}\overline{\bigcup_{t\geq s}S(t){\mathcal{B}}_1},\label{6-6}\end{eqnarray}
	where the closures are taken with respect to the weak topology of $H^{(1)}_{\delta}$, is the (maximal) universal attractor in $H^{(1)}_{\delta}$.
\end{Lemma}

\noindent {\bf Proof of Theorem \ref{th2-1}.} Combining (a), the continuity of semigroup, i.e., the existence of $C_0$-semigroup in lemmas \ref{le3-3}-\ref{le3-4}, and (b) dissipation to achieve attracting property in lemmas \ref{le4-1}-\ref{le5-2}, and (c) compactness via compact embedding in lemmas \ref{le6-1}-\ref{le6-5}, we conclude that the $\omega$-limit sets are the global attractors. This finishes the proof of Theorem \ref{th2-1}.\hspace*{\fill}$\Box$

\section{Conclusion}

This compressible micropolar fluid model is a version of the classic Navier-Stokes system coupled with an equation from microfluid models, which already leads to one of many theoretical efforts for polymeric fluids. From PDE analysis point of view, defined on the symmetric geometry (i.e., domains have spherical/cylindrical symmetry), model can be converted to a 1D system via coordinate transforms, the existence of attractor implies the long time asymptotic behavior and stability of solutions. However, since the 3D Navier-stokes equation is open, we may not expect the system \eqref{1-1} in 3D can attain better results such as existence, uniqueness, and regularity of global strong solutions than the compressible Navier-Stokes system; and the corresponding dynamic problem in domains of higher dimension is still open.\\

\section*{Acknowledgment}

This research was supported by NSFC (No. 11501199)  and the Natural Science Foundation of Henan Province (No. 18B110010). Xinguang Yang was partly supported by the Mainstay Fund from Henan Normal University. Yongjin Lu was partially supported by United States National Science Foundation (Award No. 1601127).

\end{document}